\input amstex
\input amsppt.sty
\magnification=\magstep1
\vsize=22.2truecm
\baselineskip=16truept
\nologo
\pageno=1
\topmatter
\def\Z{\Bbb Z}
\def\N{\Bbb N}

\def\bg{\bigg}
\def\({\bg(}
\def\[{\bg[}
\def\){\bg)}
\def\]{\bg]}

\def\mo{\roman{mod}}

\def\bi{\binom}
\def\Str#1#2{\left[\matrix#1\\#2\endmatrix\right]}

\def\bi{\binom}

\def\Proof{\noindent{\it Proof}}
\def\Remark{\medskip\noindent{\it  Remark}}
\def\Ack{\medskip\noindent {\bf Acknowledgment}}
\def\pmod #1{\ (\roman{mod}\ #1)}
\def\qbinom #1#2#3{\left[\matrix#1\\#2\endmatrix\right]_{#3}}
\def\fq{\roman{Q}}
\title
A $q$-analogue of Lehmer's congruence
\endtitle
\abstract
We establish the $q$-analogue of a classical congruence of Lehmer. Also, the $q$-analogues of two congruences of Morley and Granville are given.
\endabstract
\author
Hao Pan
\endauthor
\address
Department of Mathematics, Nanjing University,
Nanjing 210093, People's Republic of China
\endaddress
\email{haopan79\@yahoo.com.cn}\endemail
\endtopmatter
\document
\TagsOnRight
\heading
1. Introduction
\endheading
In 1938, Lehmer [Leh] established an interesting congruence as follows:
$$
\sum_{j=1}^{(p-1)/2}\frac{1}{j}\equiv-2\fq_p(2)+\fq_p(2)^2p\ (\mo\ p^2),\tag 1.1
$$
where $p\geq 3$ is prime and $\fq_p(2)=(2^{p-1}-1)/p$.
Lehmer's congruence can be considered as an extension
of the Wolstenholme's harmonic series congruence [W]
$$
\sum_{j=1}^{p-1}\frac{1}{j}\equiv0\ (\mo\ p^2).\tag 1.2
$$
On the other hand, the $q$-analogues of some arithmetic congruences have been investigated by several authors (e.g., see [A], [F], [C], [GZ] and [PS]).
Recently, Shi and Pan [SP] proved the following $q$-analogue of (1.2):
$$
\sum_{j=1}^{p-1}\frac{1}{[j]_q}\equiv \frac{p-1}{2}(1-q)+\frac{p^2-1}{24}(1-q)^2[p]_q\ (\mo [p]_q^2),\tag 1.3
$$
where $[n]=(1-q^n)/(1-q)=1+q+\cdots+q^{n-1}$. Obviously $(1.2)$ is deduced from $(1.3)$ when $q\to 1$.

The main purpose of the present paper is to establish the $q$-analugue of Lehmer's congruence. Set
$$
(a;q)_n=\cases (1-a)(1-aq)\cdots(1-aq^{n-1})\qquad&\text{ if } n\geq 1,\\1&\text{ if } n=0.\endcases
$$
It is easy to see that for any $m\geq 0$ with $p\nmid m$ we have a $q$-analogue of Fermat's little theorem
$$
\frac{(q^m;q^m)_{p-1}}{(q;q)_{p-1}}\equiv 1\ (\mo\ [p]_q).\tag 1.4
$$
Indeed, since
$$
[m]_q=\frac{1-q^m}{1-q}\equiv \frac{1-q^n}{1-q}=[n]_q\ (\mo\ [p]_q)
$$
whenever $m\equiv n\ (\mo\ p)$,
$$
\frac{(q^m;q^m)_{p-1}}{(q;q)_{p-1}}=\prod_{j=1}^{p-1}\frac{1-q^{jm}}{1-q^j}=\prod_{j=1}^{p-1}\frac{[jm]_q}{[j]_q}\equiv 1\ (\mo\ [p]_q).
$$
So we can define the $q$-Fermat quotient by
$$
\fq_p(m,q)=\frac{(q^m;q^m)_{p-1}/(q;q)_{p-1}-1}{[p]_q}.
$$
\proclaim{Theorem 1.1} Let $p$ be an odd prime. We have
$$
\aligned
&2\sum_{j=1}^{(p-1)/2}\frac{1}{[2j]_q}+2\fq_p(2,q)-\fq_p(2,q)^2[p]_q\\
\equiv&\(\fq_p(2,q)(1-q)+\frac{p^2-1}{8}(1-q)^2\)[p]_q\ (\mo\ [p]_q^2).
\endaligned\tag 1.5
$$
\endproclaim
In 1895, with the help of De Moivre's Theorem, Morley [M] proved that
$$
(-1)^{\frac{p-1}{2}}\bi{p-1}{(p-1)/2}\equiv 4^{p-1}\ (\mo\ p^3).\tag 1.6
$$
for any prime $p\geq 5$. In [G1], Granville generalized the congruence of Morley and showed that
$$
(-1)^{(p-1)(m-1)/2}\prod_{k=1}^{m-1}\bi{p-1}{\left\lfloor kp/m\right\rfloor}\equiv m^p-m+1\ (\mo\ p^2),\tag 1.7
$$
for any $m\geq 2$, where $\left\lfloor x\right\rfloor$ denotes the greatest integer not exceeding $x$. Now we can give the $q$-analogues of (1.6) and (1.7).
For any $m,n\in\N$, define the $q$-binomial coefficients by
$$
\qbinom{n}{m}{q}=\frac{(q;q)_n}{(q;q)_m(q;q)_{n-m}}
$$
if $n\geq m$, and if $n<m$, then we let $\qbinom{n}{m}{q}=0$.
It is easy to see that $\qbinom{n}{m}{q}$ is a polynomial in $q$ with integral coefficients, since the $q$-binomial coefficients satisfy the recurrence relation
$$
\qbinom{n+1}{m}{q}=q^m\qbinom{n}{m}{q}+\qbinom{n}{m-1}{q}.
$$
\proclaim{Theorem 1.2}
$$
(-1)^{\frac{p-1}{2}}q^{\frac{p^2-1}{4}}\Str{p-1}{(p-1)/2}_{q^2}\equiv(-q;q)_{p-1}^2-\frac{p^2-1}{24}(1-q)^2[p]_q^2\
(\mo\ [p]_q^3)\tag 1.8
$$
for any prime $p\geq 5$.
\endproclaim
\proclaim{Theorem 1.3} Let $p\geq 5$ be a prime and $m\geq 2$ be an integer with $p\nmid m$. Then
$$
(-1)^{(p-1)(m-1)/2}q^M\prod_{k=1}^{m-1}\Str{p-1}{\lfloor kp/m\rfloor}_{q^m}
\equiv\frac{m(q^m;q^m)_{p-1}}{(q;q)_{p-1}}-m+1\ (\mo\ [p]_q^2),\tag 1.9
$$
where
$$
M=m\sum_{k=1}^{m-1}\bi{\lfloor kp/m\rfloor+1}{2}.
$$
\endproclaim
The proofs of Theorems 1.1, 1.2 and 1.3 will be given in the next sections.
\heading
2. Some Lemmas
\endheading
In this section we assume that $p$ is a prime greater than $3$. And the following lemmas will be used in the proofs of Theorems 1.1 and 1.2.
\proclaim{Lemma 2.1}
$$
\sum_{j=1}^{p-1}\frac{1}{[j]_q}\equiv\frac{p-1}{2}(1-q)\pmod{[p]_q},\tag 2.1
$$
$$
\sum_{j=1}^{p-1}\frac{q^j}{[j]_q^2}\equiv-\frac{p^2-1}{12}(1-q)^2\ (\mo\ [p]_q)\tag 2.2
$$
and
$$
\sum_{j=1}^{p-1}\frac{1}{[j]_q^2}\equiv-\frac{(p-1)(p-5)}{12}(1-q)^2\ (\mo\ [p]_q).\tag 2.3
$$
\endproclaim
\Proof. See Theorem 4 in [A] and Lemma 2 in [SP].
\proclaim{Lemma 2.2}
$$
q^{kp}=\sum_{j=0}^k(-1)^j\bi{k}{j}(1-q)^j[p]_q^j.
$$
\endproclaim
\Proof.
$$
\sum_{j=0}^k(-1)^j\bi{k}{j}(1-q)^j[p]_q^j=(1-(1-q)[p]_q)^k=(1-(1-q^p))^k=q^{kp}.
$$
\qed

From Lemma 2.2, we deduce that
$$
q^{kp}\equiv1-k(1-q)[p]_q+\frac{k(k-1)}{2}(1-q)^2[p]_q^2\pmod{[p]_q^3}.\tag 2.4
$$
\proclaim{Lemma 2.3}
$$
\aligned
&4\sum_{1\leq j<k\leq p-1}\frac{(-1)^k}{[j]_q[k]_q}\\
\equiv&\(\sum_{j=1}^{p-1}\frac{(-1)^j}{[j]_q}\)^2+(p-3)(1-q)\sum_{k=1}^{p-1}\frac{(-1)^{k}}{[k]_q}+\frac{(p-1)(p+7)}{12}(1-q)^2\ (\mo\ [p]_q).
\endaligned\tag 2.5
$$
\endproclaim
\Proof. Since $p$ is odd,
$$
\aligned
\(\sum_{j=1}^{p-1}\frac{(-1)^j}{[j]_q}\)^2
=&\(\sum_{j=1}^{p-1}\frac{(-1)^j}{[j]_q}\)\(\sum_{j=1}^{p-1}\frac{(-1)^j}{[j]_q}-(1-q)\sum_{j=1}^{p-1}(-1)^j\)\\
=&\(\sum_{j=1}^{p-1}\frac{(-1)^j}{[j]_q}\)\(\sum_{j=1}^{p-1}\frac{(-q)^j}{[j]_q}\)\\
=&\sum_{k=2}^{2p-2}(-1)^k\sum_{j=\max\{1,k-p+1\}}^{\min\{k-1,p-1\}}\frac{q^j}{[j]_q[k-j]_q}.
\endaligned
$$
Then we have
$$
\aligned
&\(\sum_{j=1}^{p-1}\frac{(-1)^j}{[j]_q}\)\(\sum_{j=1}^{p-1}\frac{(-q)^j}{[j]_q}\)-(-1)^p\sum_{j=1}^{p-1}\frac{q^j}{[j]_q[p-j]_q}\\
=&\sum_{k=2}^{p-1}(-1)^k\sum_{j=1}^{k-1}\frac{q^j}{[j]_q[k-j]_q}+\sum_{k=p+1}^{2p-2}(-1)^k\sum_{j=k-p+1}^{p-1}\frac{q^j}{[j]_q[k-j]_q}\\
=&\sum_{k=2}^{p-1}(-1)^k\sum_{j=1}^{k-1}\frac{q^j}{[j]_q[k-j]_q}+\sum_{l=2}^{p-1}(-1)^{2p-l}\sum_{j=p-l+1}^{p-1}\frac{q^j}{[j]_q[2p-l-j]_q}\\
&(\text{here } l=2p-k)\\
=&\sum_{k=2}^{p-1}(-1)^k\sum_{j=1}^{k-1}\frac{q^j}{[j]_q[k-j]_q}+\sum_{l=2}^{p-1}(-1)^{l}\sum_{i=1}^{l-1}\frac{q^{p+i-l}}{[p+i-l]_q[p-i]_q}.\\
&(\text{here } i=l+j-p)
\endaligned
$$
Note that
$$
\frac{q^{p+i-l}}{[p+i-l]_q[p-i]_q}=\frac{q^{p+i-l}(1-q)^2}{(1-q^{p+i-l})(1-q^{p-i})}\equiv\frac{q^i(1-q)^2}{(1-q^{l-i})(1-q^i)}\ (\mo\ [p]_q).
$$
It follows that
$$
\aligned
&\(\sum_{j=1}^{p-1}\frac{(-1)^j}{[j]_q}\)\(\sum_{j=1}^{p-1}\frac{(-q)^j}{[j]_q}\)+\sum_{j=1}^{p-1}\frac{q^j}{[j]_q[p-j]_q}\\
\equiv&\sum_{k=2}^{p-1}(-1)^k\sum_{j=1}^{k-1}\frac{q^j}{[j]_q[k-j]_q}+\sum_{l=2}^{p-1}(-1)^{l}\sum_{i=1}^{l-1}\frac{q^i}{[i]_q[l-i]_q}\\
=&2\sum_{k=2}^{p-1}(-1)^k\sum_{j=1}^{k-1}\frac{q^j(1-q)^2}{(1-q^j)(1-q^{k-j})}\\
=&2\sum_{k=2}^{p-1}\frac{(-1)^kq^k(1-q)^2}{1-q^k}\sum_{j=1}^{k-1}\(\frac{1}{q^{k-j}-q^k}+\frac{1}{1-q^{k-j}}\)\\
=&2\sum_{1\leq j<k\leq p-1}\frac{(-1)^k(q^k+q^j)}{[j]_q[k]_q}\ (\mo\ [p]_q).
\endaligned
$$
We can write
$$
\aligned
&\sum_{1\leq j<k\leq p-1}\frac{(-1)^k(q^k+q^j)}{[j]_q[k]_q}\\
=&\sum_{1\leq j<k\leq p-1}\frac{(-1)^k(2-(1-q^k)-(1-q^j))}{[j]_q[k]_q}\\
=&2\sum_{1\leq j<k\leq p-1}\frac{(-1)^k}{[j]_q[k]_q}-(1-q)\(\sum_{k=2}^{p-1}\frac{(-1)^k(k-1)}{[k]_q}+\sum_{j=1}^{p-2}\frac{1}{[j]_q}\sum_{k=j+1}^{p-1}(-1)^k\).
\endaligned
$$
Now
$$
\aligned
\sum_{k=2}^{p-1}\frac{(-1)^k(k-1)}{[k]_q}=&\frac{1}{2}\sum_{k=1}^{p-1}\(\frac{(-1)^k(k-1)}{[k]_q}+\frac{(-1)^{p-k}(p-k-1)}{[p-k]_q}\)\\
=&\frac{1}{2}\sum_{k=1}^{p-1}(-1)^k(k-1)\(\frac{1}{[k]_q}+\frac{1}{[p-k]_q}\)+\frac{p-2}{2}\sum_{k=1}^{p-1}\frac{(-1)^{k}}{[k]_q}\\
=&\frac{1}{2}\sum_{k=1}^{p-1}(-1)^k(k-1)\(\frac{[p]_q}{[k][p-k]_q}+(1-q)\)+\frac{p-2}{2}\sum_{k=1}^{p-1}\frac{(-1)^{k}}{[k]_q}\\
\equiv&\frac{p-1}{4}(1-q)+\frac{p-2}{2}\sum_{k=1}^{p-1}\frac{(-1)^{k}}{[k]_q}\ (\mo\ [p]_q).
\endaligned
$$
And from (2.1) we have
$$
\aligned
\sum_{j=1}^{p-2}\frac{1}{[j]_q}\sum_{k=j+1}^{p-1}(-1)^k=&\frac{1}{2}\sum_{j=1}^{p-1}\frac{1-(-1)^j}{[j]_q}\\
\equiv&\frac{(p-1)}{4}(1-q)-\frac{1}{2}\sum_{j=1}^{p-1}\frac{(-1)^j}{[j]_q}\ (\mo\ [p]_q).
\endaligned
$$
Finally by (2.3),
$$
\aligned
\sum_{j=1}^{p-1}\frac{q^j}{[j]_q[p-j]_q}=\sum_{j=1}^{p-1}\frac{q^{p-j}}{[j]_q[p-j]_q}=&\sum_{j=1}^{p-1}\frac{q^p}{[j]_q([p]-[j]_q)}\\
\equiv&\frac{(p-1)(p-5)}{12}(1-q)^2\ (\mo\ [p]_q).
\endaligned
$$
Thus combining the equations and congruences above, we obtain that
$$
\aligned
&4\sum_{1\leq j<k\leq p-1}(-1)^k\frac{1}{[j]_q[k]_q}-\(\sum_{j=1}^{p-1}\frac{(-1)^j}{[j]_q}\)\(\sum_{j=1}^{p-1}\frac{(-q)^j}{[j]_q}\)\\
\equiv&(p-3)(1-q)\sum_{k=1}^{p-1}\frac{(-1)^{k}}{[k]_q}+\frac{(p-1)(p+7)}{12}(1-q)^2\ (\mo\ [p]_q).
\endaligned
$$
\qed
\proclaim{Lemma 2.4}
$$
\aligned
\sum_{j=1}^{p-1}\frac{(-1)^j}{[j]_q}\equiv 2\sum_{j=1}^{(p-1)/2}\frac{1}{[2j]_q}-\frac{p-1}{2}(1-q)-\frac{p^2-1}{24}(1-q)^2[p]_q\ (\mo\ [p]_q^2).
\endaligned\tag 2.6
$$
\endproclaim
\Proof. Clearly
$$
\sum_{j=1}^{p-1}\frac{(-1)^j}{[j]_q}=\sum_{j=1}^{(p-1)/2}\frac{1}{[2j]_q}-\sum_{j=1}^{(p-1)/2}\frac{1}{[2j-1]_q}
=\sum_{j=1}^{(p-1)/2}\frac{1}{[2j]_q}-\sum_{j=1}^{(p-1)/2}\frac{1}{[p-2j]_q}.
$$
Observe that
$$
\frac{1}{[p-2j]_q}=\frac{q^{2j}}{[p]_q-[2j]_q}=\frac{q^{2j}([p]_q+[2j]_q)}{[p]_q^2-[2j]_q^2}\equiv-\frac{q^{2j}([p]_q+[2j]_q)}{[2j]_q^2}\ (\mo\ [p]_q^2).
$$
And by (2.2), we have
$$
\aligned
-\frac{p^2-1}{12}(1-q)^2\equiv
\sum_{j=1}^{p-1}\frac{q^j}{[j]_q^2}
=&\sum_{j=1}^{(p-1)/2}\frac{q^{2j}}{[2j]_q^2}+\sum_{j=1}^{(p-1)/2}\frac{q^{p-2j}}{[p-2j]_q^2}\\
=&\sum_{j=1}^{(p-1)/2}\frac{q^{2j}}{[2j]_q^2}+\sum_{j=1}^{(p-1)/2}\frac{q^{p+2j}}{([p]_q-[2j]_q)^2}\\
\equiv&2\sum_{j=1}^{(p-1)/2}\frac{q^{2j}}{[2j]_q^2}\ (\mo\ [p]_q).
\endaligned\tag 2.7
$$
Hence
$$
\aligned
\sum_{j=1}^{p-1}\frac{(-1)^j}{[j]_q}\equiv&\sum_{j=1}^{(p-1)/2}\frac{1}{[2j]_q}+\sum_{j=1}^{(p-1)/2}\frac{q^{2j}([p]_q+[2j]_q)}{[2j]_q^2}\\
\equiv&\sum_{j=1}^{(p-1)/2}\frac{1+q^{2j}}{[2j]_q}-\frac{p^2-1}{24}(1-q)^2[p]_q\\
=&2\sum_{j=1}^{(p-1)/2}\frac{1}{[2j]_q}-\frac{p-1}{2}(1-q)-\frac{p^2-1}{24}(1-q)^2[p]_q\ (\mo\ [p]_q^2).
\endaligned
$$
\qed
\heading 3. Proofs of Theorems 1.1 and 1.2
\endheading
\medskip\noindent{\it  Proof of Theorem 1.1}. One can directly verify (1.5) when $p=3$. So below we assume that $p\geq 5$.
It is well-known (cf. Corollary 10.2.2 of [AAR]) that
$$
(x;q)_n=\sum_{j=0}^n\Str{n}{j}_qq^{\bi{j}{2}}(-x)^j.
$$
Then we have
$$
\aligned
\frac{(-1;q)_q-q^{\bi{p}{2}}-1}{[p]_q}
=&\frac{1}{[p]_q}\sum_{k=1}^{p-1}\Str{p}{k}_qq^{\bi{k}{2}}\\
=&\sum_{k=1}^{p-1}\frac{1}{[k]_q}\prod_{j=1}^{k-1}\frac{q^j(1-q^{p-j})}{1-q^j}\\
=&\sum_{k=1}^{p-1}\frac{1}{[k]_q}\prod_{j=1}^{k-1}\(\frac{[p]_q}{[j]_q}-1\)\\
\equiv&[p]_q\sum_{1\leq j<k\leq p-1}\frac{(-1)^k}{[j]_q[k]_q}-\sum_{k=1}^{p-1}\frac{(-1)^k}{[k]_q}\pmod{[p]_q^2}.
\endaligned\tag 3.1
$$
Consequently
$$
\sum_{k=1}^{p-1}\frac{(-1)^k}{[k]_q}\equiv-\frac{(-1;q)_q-q^{\bi{p}{2}}-1}{[p]_q}
\equiv-\frac{2(-q;q)_{q-1}-2}{[p]_q}-\frac{p-1}{2}(1-q)\pmod{[p]_q}.
$$
Thus applying Lemma 2.3, we have
$$
\aligned
&\sum_{1\leq j<k\leq p-1}\frac{(-1)^k}{[j]_q[k]_q}-\frac{(p-1)(p+7)}{48}(1-q)^2\\
\equiv&\frac{1}{4}\(\sum_{k=1}^{p-1}\frac{(-1)^k}{[k]_q}\)\(\sum_{k=1}^{p-1}\frac{(-1)^k}{[k]_q}+(p-3)(1-q)\)\\
\equiv&\frac{1}{4}\(-2\fq_p(2,q)-\frac{p-1}{2}(1-q)\)\(-2\fq_p(2,q)+\frac{p-5}{2}(1-q)\)\\
=&\fq_p(2,q)^2+\fq_p(2,q)(1-q)-\frac{(p-1)(p-5)}{16}(1-q)^2\pmod{[p]_q}.
\endaligned\tag 3.2
$$
On the other hand, it follows from (2.4) that
$$
\aligned
&\frac{(-1;q)_q-q^{\bi{p}{2}}-1}{[p]_q}\\
\equiv&\frac{2(-q;q)_{q-1}-2}{[p]_q}+\frac{p-1}{2}(1-q)-\frac{(p-1)(p-3)}{8}(1-q)^2[p]_q\pmod{[p]_q^2}.
\endaligned
$$
Then by Lemma 2.4,
$$
\aligned
&\sum_{k=1}^{p-1}\frac{(-1)^k}{[k]_q}+\frac{(-1;q)_q-q^{\bi{p}{2}}-1}{[p]_q}\\
\equiv&2\sum_{j=1}^{(p-1)/2}\frac{1}{[2j]_q}+2\fq_p(2,q)-\frac{(p-1)(p-2)}{6}(1-q)^2[p]_q\pmod{[p]_q^2}.
\endaligned\tag 3.3
$$
Combining (3.1), (3.2) and (3.3), the desired (1.5) is obtained.
\qed
\medskip\noindent{\it  Proof of Theorem 1.2}. Since
$$
\Str{p-1}{(p-1)/2}_{q^2}=\prod_{j=1}^{(p-1)/2}\frac{[p-j]_{q^2}}{[j]_{q^2}}=\prod_{j=1}^{(p-1)/2}\frac{[p]_{q^2}-[j]_{q^2}}{q^{2j}[j]_{q^2}},
$$
we have
$$
\aligned
&(-1)^{\frac{p-1}{2}}q^{\frac{p^2-1}{4}}\Str{p-1}{(p-1)/2}_{q^2}\\
=&\prod_{j=1}^{(p-1)/2}\(1-\frac{[p]_{q^2}}{[j]_{q^2}}\)\\
\equiv&1-\frac{1+q^p}{1+q}\sum_{j=1}^{(p-1)/2}\frac{[p]_{q}}{[j]_{q^2}}+\frac{(1+q^p)^2}{(1+q)^2}\sum_{1\leq j<k\leq (p-1)/2}\frac{[p]_{q}^2}{[j]_{q^2}[k]_{q^2}}\pmod{[p]_{q^2}^3}.\\
\endaligned\tag 3.4
$$
From Theorem 1.1, we deduce that
$$
\aligned
&\frac{1+q^p}{1+q}\sum_{j=1}^{(p-1)/2}\frac{1}{[j]_{q^2}}\\
=&(1+q^p)\sum_{j=1}^{(p-1)/2}\frac{1}{[2j]_{q}}\\
\equiv&-(1+q^p)\fq_p(2,q)+\frac{1+q^p}{2}\(\fq_p(2,q)^2+\fq_p(2,q)(1-q)+\frac{p^2-1}{8}(1-q)^2\)[p]_q\\
\equiv&-(1+q^p)\fq_p(2,q)+\fq_p(2,q)^2[p]_q+\fq_p(2,q)(1-q)[p]_q
+\frac{p^2-1}{8}(1-q)^2[p]_q\\
&\pmod{[p]_q^2}.
\endaligned\tag 3.5
$$
Notice that
$$
\aligned
\sum_{1\leq j<k\leq (p-1)/2}\frac{1}{[j]_{q^2}[k]_{q^2}}=&\frac{1}{2}\(\(\sum_{j=1}^{(p-1)/2}\frac{1}{[j]_{q^2}}\)^2-\sum_{j=1}^{(p-1)/2}\frac{1}{[j]_{q^2}^2}\)\\
=&\frac{(1+q)^2}{2}\(\(\sum_{j=1}^{(p-1)/2}\frac{1}{[2j]_{q}}\)^2-\sum_{j=1}^{(p-1)/2}\frac{1}{[2j]_{q}^2}\).
\endaligned
$$
And Theorem 1.1 implies that
$$
\sum_{j=1}^{(p-1)/2}\frac{1}{[2j]_q}\equiv-\fq_p(2,q)\pmod{[p]_q}.
$$
Then using (2.7),
$$
\aligned
\sum_{j=1}^{(p-1)/2}\frac{1}{[2j]_{q}^2}
=&\sum_{j=1}^{(p-1)/2}\frac{q^{2j}}{[2j]_{q}^2}+(1-q)\sum_{j=1}^{(p-1)/2}\frac{1}{[2j]_{q}}\\
\equiv&-\frac{p^2-1}{24}(1-q)^2-\fq_p(2,q)(1-q)\pmod{[p]_q}.
\endaligned
$$
Consequently
$$
\aligned
&\frac{2}{(1+q)^2}\sum_{1\leq j<k\leq (p-1)/2}\frac{1}{[j]_{q^2}[k]_{q^2}}\\
\equiv&\fq_p(2,q)^2+\fq_p(2,q)(1-q)+\frac{p^2-1}{24}(1-q)^2\pmod{[p]_q}.\\
\endaligned\tag 3.6
$$
Thus it follows from (3.4), (3.5) and (3.6) that
$$
\aligned
&(-1)^{\frac{p-1}{2}}q^{\frac{p^2-1}{4}}\Str{p-1}{(p-1)/2}_{q^2}-1\\
\equiv&[p]_{q}^2\cdot\frac{4}{(1+q)^2}\sum_{1\leq j<k\leq (p-1)/2}\frac{1}{[j]_{q^2}[k]_{q^2}}-[p]_{q}\cdot\frac{1+q^p}{1+q}\sum_{j=1}^{(p-1)/2}\frac{1}{[j]_{q^2}}\\
\equiv&\fq_p(2,q)^2[p]_q^2+(1+q^p)\fq_p(2,q)[p]_q+\fq_p(2,q)(1-q)[p]_q^2-\frac{p^2-1}{24}(1-q)^2[p]_q^2\\
=&((-q;q)_{p-1}-1)((-q;q)_{p-1}+1)-\frac{p^2-1}{24}(1-q)^2[p]_q^2\pmod{[p]_q^3}.
\endaligned
$$
\qed
\heading
4. Fermat Quotient
\endheading
\proclaim{Lemma 4.1} Let $p$ be an odd prime. Suppose that $m$ is a positive integer with $(m,p)=1$. Then
$$
\fq_p(m,q)\equiv\sum_{j=1}^{p-1}\frac{\lfloor jm/p\rfloor}{[jm]_q}-\frac{(p-1)(m-1)}{2}(1-q)\ (\mo\ [p]_q).\tag 4.1
$$
\endproclaim
\Proof. For each $j\in\{1,2,\ldots,p-1\}$, let
$$
r_j=jm-\left\lfloor jm/p\right\rfloor p.
$$
Then
$$
\aligned
\frac{(q^m;q^m)_{p-1}}{(q;q)_{p-1}}=\prod_{j=1}^{p-1}\frac{1-q^{jm}}{1-q^j}
=&\prod_{j=1}^{p-1}\(\frac{1-q^{r_j}}{1-q^j}+\frac{q^{r_j}(1-q^{\lfloor jm/p\rfloor p})}{1-q^j}\)\\
=&\prod_{j=1}^{p-1}\frac{1-q^{r_j}}{1-q^j}\(1+\frac{q^{r_j}(1-q^{\lfloor jm/p\rfloor p})}{1-q^{r_j}}\)
\endaligned
$$
Since $r_j$ runs through $1,2,\ldots,p-1$ as $j$ does so, we have
$$
\aligned
\frac{(q^m;q^m)_{p-1}}{(q;q)_{p-1}}=&\prod_{j=1}^{p-1}\(1+\frac{q^{r_j}(1-q^{\lfloor jm/p\rfloor p})}{1-q^{r_j}}\)\\
\equiv&1+(1-q^p)\sum_{j=1}^{p-1}\frac{q^{r_j}}{1-q^{r_j}}\cdot\frac{1-q^{\lfloor jm/p\rfloor p}}{1-q^p}\\
\equiv&1+(1-q^p)\sum_{j=1}^{p-1}\left\lfloor \frac{jm}{p}\right\rfloor\frac{q^{r_j}}{1-q^{r_j}}\\
\equiv&1+[p]_q\sum_{j=1}^{p-1}\left\lfloor \frac{jm}{p}\right\rfloor\frac{q^{jm}}{[jm]_q}\ (\mo\ [p]_q^2).
\endaligned
$$
Finally,
$$
\aligned
\sum_{j=1}^{p-1}\left\lfloor \frac{jm}{p}\right\rfloor\frac{q^{jm}}{[jm]_q}=&\sum_{j=1}^{p-1}\frac{\lfloor jm/p\rfloor}{[jm]_q}-(1-q)\sum_{j=1}^{p-1}\left\lfloor \frac{jm}{p}\right\rfloor\\
=&\sum_{j=1}^{p-1}\frac{\lfloor jm/p\rfloor}{[jm]_q}-\frac{(p-1)(m-1)}{2}(1-q).
\endaligned
$$
We are done.\qed

\Remark. Letting $q\to1$ in (4.1), we obtain that
$$
\frac{m^p-m}{p}\equiv\sum_{j=1}^{p-1}\frac{\lfloor jm/p\rfloor}{j}\pmod{p},
$$
which was firstly discovered by Lerch [Ler].
\medskip\noindent{\it  Proof of Theorem 1.3}. We write
$$
\Str{p-1}{\lfloor kp/m\rfloor}_{q^m}=\prod_{j=1}^{\lfloor kp/m\rfloor}\frac{[p]_{q^m}-[j]_{q^m}}{q^{jm}[j]_{q^m}}=q^{-m\bi{\lfloor kp/m\rfloor+1}{2}}\prod_{j=1}^{\lfloor kp/m\rfloor}\(\frac{[p]_{q^m}}{[j]_{q^m}}-1\).
$$
As $p\nmid m$, $[p]_q$ divides $[p]_{q^m}=(1-q^{mp})/(1-q^m)$. Thus
$$
\aligned
&(-1)^{(p-1)(m-1)/2}q^{\sum_{k=1}^{m-1}m\bi{\lfloor kp/m\rfloor+1}{2}}\prod_{k=1}^{m-1}\Str{p-1}{\lfloor kp/m\rfloor}_{q^m}\\
=&\prod_{k=1}^{m-1}\prod_{j=1}^{\lfloor kp/m\rfloor}\(1-\frac{[p]_{q^m}}{[j]_{q^m}}\)\\
\equiv&1-[p]_{q^m}\sum_{k=1}^{m-1}\sum_{1\leq j<kp/m}\frac{1}{[j]_{q^m}}\\
&(\text{here }kp/m\not\in\Z,\ \text{so }j\leq\lfloor kp/m\rfloor<kp/m)\\
=&1-[p]_{q^m}\sum_{j=1}^{p-1}\frac{m-1-\lfloor jm/p\rfloor}{[j]_{q^m}}\ (\mo\ [p]_{q}^2).
\endaligned
$$
In view of (2.1) and Lemma 4.1, we have
$$
\aligned
&[p]_{q^m}\sum_{j=1}^{p-1}\frac{m-1-\lfloor jm/p\rfloor}{[j]_{q^m}}\\
=&(m-1)[mp]_{q}\sum_{j=1}^{p-1}\frac{1}{[jm]_{q}}-[mp]_{q}\sum_{j=1}^{p-1}\frac{\lfloor jm/p\rfloor}{[jm]_q}\\
\equiv&(m-1)[mp]_{q}\cdot\frac{(p-1)}{2}(1-q)-[mp]_q\fq_p(m,q)-\frac{(p-1)(m-1)}{2}(1-q)[mp]_q\\
\equiv&-m[p]_q\fq_p(m,q)\ (\mo\ [p]_q^2).
\endaligned
$$
This concludes our proof.\qed

\Remark. For further developments of Granville's congruence (1.7), the reader is referred to [S].

\heading
5. A conjecture of Skula
\endheading
Recently with help of polynomials over finite fields, Granville [G2] confirmed a conjecture of Skula:
$$
\(\frac{2^{p-1}-1}{p}\)^2\equiv-\sum_{j=1}^{p-1}\frac{2^j}{j^2}\pmod{[p]_q}\tag 5.1
$$
for any prime $p\geq 5$. Using our $q$-analogue of Lehmer's congruence, we have the following $q$-analogue of (5.1):
\proclaim{Theorem 5.1} Let $p\geq 5$ be a prime. Then
$$
\aligned
&\sum_{j=1}^{p-1}\frac{q^j(-q;q)_j}{[j]_q^2}+\fq_p(2,q)^2\\
\equiv&-(p-1)\fq_p(2,q)(1-q)-\frac{(7p-5)(p-1)}{24}(1-q)^2\pmod{[p]_q}.
\endaligned\tag 5.2
$$
\endproclaim
\proclaim{Lemma 5.2}
$$
\sum_{k=0}^n(-1)^k\qbinom{n}{k}{q}q^{\bi{n-k}{2}}(-q;q)_k=(-1)^nq^{\bi{n+1}{2}}.\tag 5.3
$$
\endproclaim
\Proof. From the well-known $q$-binomial theorem (cf. Theorem 10.2.1 of [AAR]), we have
$$
\sum_{k=0}^\infty\frac{(-1)^kq^{\bi{k}{2}}}{(q;q)_k}x^k=(x;q)_{\infty}
$$
and
$$
\sum_{k=0}^\infty\frac{(-q;q)_k}{(q;q)_k}x^k=\frac{(-qx;q)_{\infty}}{(x;q)_{\infty}}.
$$
Then by comparing the coefficient of $x^n$ in the both sides of
$$
(x;q)_{\infty}\cdot\frac{(-qx;q)_{\infty}}{(x;q)_{\infty}}=(-qx;q)_{\infty},
$$
we obtain that
$$
\sum_{k=0}^n\frac{(-1)^{n-k}q^{\bi{n-k}{2}}(-q;q)_k}{(q;q)_{n-k}(q;q)_k}=\frac{q^{\bi{n}{2}+n}}{(q;q)_n},
$$
which is an equivalent form of (5.3).
\qed
\proclaim{Corollary 5.3} For any odd prime $p$, we have
$$
\sum_{j=1}^{p-1}\frac{q^j(-q;q)_j}{[j]_q}\equiv-2\fq_p(2,q)-(p-1)(1-q)\pmod{[p]_q}.\tag 5.4
$$
\endproclaim
\Proof. From Lemma 5.2, we deduce that
$$
\aligned
\sum_{j=1}^{p-1}\frac{q^j(-q;q)_j}{[j]_q}\equiv&\sum_{j=1}^{p-1}\frac{q^{p(p-1)/2-jp+j}(-q;q)_j}{[j]_q}\\
\equiv&-\frac{1}{[p]_q}\sum_{j=1}^{p-1}(-1)^jq^{\bi{p}{2}+\bi{j}{2}-jp+j}\qbinom{p}{j}{q}(-q;q)_j\\
=&-\frac{1}{[p]_q}\sum_{j=1}^{p-1}(-1)^jq^{\bi{p-j}{2}}\qbinom{p}{j}{q}(-q;q)_j\\
=&-\frac{1}{[p]_q}((-1)^pq^{\bi{p+1}{2}}-q^{\bi{p}{2}}-(-1)^p(-q;q)_p)\pmod{[p]_q}.
\endaligned
$$
Notice that
$$
\frac{2-q^{\bi{p+1}{2}}-q^{\bi{p}{2}}}{[p]_q}\equiv \frac{p+1}{2}(1-q)+\frac{p-1}{2}(1-q)=p(1-q)\pmod{[p]_q},
$$
and that
$$
\aligned
\frac{(-q;q)_p-2}{[p]_q}=\frac{(-q;q)_{p-1}(1+q^p)-2}{[p]_q}=&\frac{2(-q;q)_{p-1}-2}{[p]_q}-(1-q)(-q;q)_{p-1}\\
\equiv&2\fq_p(2,q)-(1-q)\pmod{[p]_q}.
\endaligned
$$
Hence
$$
\aligned
\sum_{j=1}^{p-1}\frac{q^j(-q;q)_j}{[j]_q}\equiv&-\frac{(-q;q)_p-q^\bi{p+1}{2}-q^{\bi{p}{2}}}{[p]_q}\\
\equiv&-2\fq_p(2,q)-(p-1)(1-q)\pmod{[p]_q}.
\endaligned
$$
We are done.
\qed
\Remark. Corollary 5.3 is the $q$-analogue of an observation of Glashier:
$$
\frac{2^{p-1}-1}{p}\equiv -\sum_{j=1}^{p-1}\frac{2^{j-1}}{j}\pmod{p}.\tag 5.5
$$
\proclaim{Lemma 5.4}
$$
\sum_{k=1}^n(-1)^k\qbinom{n}{k}{q}q^{\bi{n-k}{2}}\frac{(-q;q)_k}{[k]_q}=q^{\bi{n}{2}}\sum_{k=1}^n\frac{(-q)^k-1}{[k]_q}.\tag 5.6
$$
\endproclaim
\Proof.
We make an induction on $n$. The case $n=1$ is trivial. Assume that $n>1$ and that (5.6) holds for the smaller values of $n$. Then we conclude that
$$
\aligned
&\sum_{k=1}^n(-1)^k\qbinom{n}{k}{q}q^{\bi{n-k}{2}}\frac{(-q;q)_k}{[k]_q}\\
=&\sum_{k=1}^n(-1)^k\(q^k\qbinom{n-1}{k}{q}+\qbinom{n-1}{k-1}{q}\)q^{\bi{n-k}{2}}\frac{(-q;q)_k}{[k]_q}\\
=&q^{n-1}\sum_{k=1}^{n-1}(-1)^k\qbinom{n-1}{k}{q}q^{\bi{n-k-1}{2}}\frac{(-q;q)_k}{[k]_q}+\frac{1}{[n]_q}\sum_{k=1}^n(-1)^k\qbinom{n}{k}{q}q^{\bi{n-k}{2}}(-q;q)_k\\
=&q^{\bi{n}{2}}\sum_{k=1}^{n-1}\frac{(-q)^k-1}{[k]_q}+\frac{1}{[n]_q}((-1)^nq^{\bi{n+1}{2}}-q^{\bi{n}{2}}),
\endaligned
$$
where in the last step we apply the induction hypothesis and Lemma 5.2.
\qed

\medskip\noindent{\it  Proof of Theorem 5.1}. Using Lemma 5.4, we have
$$
\aligned
\sum_{j=1}^{p-1}\frac{q^j(-q;q)_j}{[j]_q^2}\equiv&\sum_{j=1}^{p-1}\frac{q^{p(p-1)/2-jp+j}(-q;q)_j}{[j]_q^2}\\
\equiv&-\frac{1}{[p]_q}\sum_{j=1}^{p-1}(-1)^jq^{\bi{p-j}{2}}\qbinom{p}{j}{q}\frac{(-q;q)_j}{[j]_q}\\
=&-\frac{q^{\bi{p}{2}}}{[p]_q}\sum_{j=1}^p\frac{(-q)^j-1}{[j]_q}+(-1)^p\frac{(-q;q)_p}{[p]_q^2}\\
=&-\frac{q^{\bi{p}{2}}}{[p]_q}\sum_{j=1}^{p-1}\frac{(-q)^j-1}{[j]_q}-\frac{(-q;q)_p-q^{\bi{p+1}{2}}-q^{\bi{p}{2}}}{[p]_q^2}\pmod{[p]_q}.
\endaligned
$$
With help of (1.3) and Theorem 1.1,
$$
\aligned
\sum_{k=1}^{p-1}\frac{(-q)^j-1}{[j]_q}=&-\sum_{k=1}^{p-1}\frac{(-1)^j(1-q^{j})}{[j]_q}+\sum_{k=1}^{p-1}\frac{(-1)^j-1}{[j]_q}\\
=&-\sum_{j=1}^{(p-1)/2}\frac{2}{[2j-1]_q}\\
=&\sum_{j=1}^{(p-1)/2}\frac{2}{[2j]_q}-\sum_{k=1}^{p-1}\frac{2}{[j]_q}\\
\equiv&-2\fq_p(2,q)+\fq_p(2,q)^2[p]_q+\fq_p(2,q)(1-q)[p]_q\\
&+\frac{p^2-1}{8}(1-q)^2[p]_q-((p-1)(1-q)+\frac{p^2-1}{12}(1-q)^2[p]_q)\\
&\pmod{[p]_q^2}.
\endaligned
$$
And by (2.4) we have
$$
\aligned
&\frac{(-q;q)_p-q^{\bi{p+1}{2}}-q^{\bi{p}{2}}}{[p]_q^2}\\
=&\frac{2(-q;q)_{p-1}-q^{\bi{p+1}{2}}-q^{\bi{p}{2}}}{[p]_q^2}-\frac{(-q;q)_{p-1}}{[p]_q}(1-q)\\
\equiv&\frac{2(-q;q)_{p-1}-2}{[p]_q^2}+\frac{p}{[p]_q}(1-q)-\frac{(p-1)^2}{4}(1-q)^2-\frac{(-q;q)_{p-1}}{[p]_q}(1-q)\\
=&\frac{2\fq_p(2,q)}{[p]_q}+\frac{(p-1)(1-q)}{[p]_q}-\frac{(p-1)^2}{4}(1-q)^2-\fq_p(2,q)(1-q)\pmod{[p]_q}.
\endaligned
$$
Therefore
$$
\aligned
&\sum_{j=1}^{p-1}\frac{q^j(-q;q)_j}{[j]_q^2}\\
\equiv&\(\frac{(p-1)(1-q)}{2}-\frac{1}{[p]_q}\)\sum_{j=1}^{p-1}\frac{(-q)^j-1}{[j]_q}-\frac{(-q;q)_p-q^{\bi{p+1}{2}}-q^{\bi{p}{2}}}{[p]_q^2}\\
\equiv&-(p-1)\fq_p(2,q)(1-q)-\fq_p(2,q)^2-\frac{(7p-5)(p-1)}{24}(1-q)^2\pmod{[p]_q}.
\endaligned
$$
\qed

\Ack. I thank my advisor, Prof. Zhi-Wei Sun, for his useful comments on this paper.

\enddocument